\documentclass[11pt,letterpaper]{article}

\usepackage[utf8]{inputenc}
\usepackage[margin=1in]{geometry}
\usepackage{amsmath,amsthm,amssymb}
\usepackage{mathtools}
\usepackage{mathrsfs}
\usepackage{xcolor}
\usepackage[unicode=true,pdfusetitle,
bookmarks=true,bookmarksnumbered=false,bookmarksopen=false,
breaklinks=false,pdfborder={0 0 0},pdfborderstyle={},backref=false,colorlinks=false]
{hyperref}

\usepackage[nameinlink,capitalise,noabbrev]{cleveref}
\usepackage{comment}

\theoremstyle{definition}
\newtheorem{definition}{Definition}[section]
\newtheorem*{remarkstar}{Remark}
\newtheorem*{definitionstar}{Definition}

\theoremstyle{plain}
\newtheorem{theorem}[definition]{Theorem}

\newtheorem{lemma}[definition]{Lemma}

\newtheorem{proposition}[definition]{Proposition}

\newtheorem{observation}[definition]{Observation}
\newtheorem{claim}[definition]{Claim}

\newtheorem{cor}[definition]{Corollary}

\renewcommand{\le}{\leqslant}
\renewcommand{\ge}{\geqslant}
\renewcommand{\leq}{\leqslant}
\renewcommand{\geq}{\geqslant}

\def \eps {\varepsilon}
\def \etan {\delta}
\def \es {\emptyset}
\def \tri {\triangle}
\def \sm {\setminus}

\def \mF {\mathcal{F}}
\def \mG {\mathcal{G}}
\def \mH {\mathcal{H}}
\def \mS {\mathcal{S}}
\def \mX {\mathcal{X}}
\def \mY {\mathcal{Y}}
\def \mU {\mathcal{U}}
\def \mV {\mathcal{V}}
\def \ce {\coloneqq}
\global\long\def\disj{\operatorname{dp}}
\global\long\def\bN{\mathbb{N}}
\newcommand{\card}[1]{\left| #1 \right|}

\title{Maximal 3-wise Intersecting Families}
\author{
J\'ozsef Balogh\thanks{Department of Mathematics, University of Illinois at Urbana-Champaign, Urbana, Illinois 61801, USA. E-mail: \texttt{\{jobal, cechen4, haoranl8\}@illinois.edu}.}  \thanks{Research is partially supported by NSF Grant DMS-1764123, NSF RTG grant DMS 1937241, Arnold O. Beckman Research Award (UIUC Campus Research Board RB 22000), and the Langan Scholar Fund (UIUC). %Haoran Luo is partially supported by UIUC Campus Research Board RB 22000.
}
\and
Ce Chen\footnotemark[1]
\and
Kevin Hendrey\thanks{Discrete Mathematics Group, Institute for Basic Science (IBS), Daejeon, Republic of Korea. E-mail: {\nolinkurl{{kevinhendrey, benlund}@ibs.re.kr}}. This work was supported by the Institute for Basic Science (IBS-R029-C1).}
\and
Ben Lund\footnotemark[3]
\and
Haoran Luo\footnotemark[1] \thanks{Research is partially supported by UIUC Campus Research Board RB 22000.}
\and
Casey Tompkins\thanks{Alfr\'ed R\'enyi Institute of Mathematics, Hungarian Academy of Sciences. E-mail: {\tt casey@renyi.hu}. This work was supported by NKFIH grant K135800.}
\and
Tuan Tran\thanks{School of Mathematical Sciences, University of Science and Technology of China, China. E-mail: {\tt trantuan@ustc.edu.cn}. This work was supported by the Institute for Basic Science (IBS-R029-Y1), and the Outstanding Young Talents Program (Overseas) of the National Natural Science Foundation of China.}
}
\date{}

\begin{document}
\maketitle

\begin{abstract}
    A family $\mF$ on ground set $[n]:=\{1,2,\ldots, n\}$ 
    is \emph{maximal $k$-wise intersecting} if every collection of at most $k$ sets in $\mF$ has non-empty intersection, and no other set can be added to $\mF$ while maintaining this property. In 1974, Erd\H{o}s and Kleitman asked for the minimum size of a maximal $k$-wise intersecting family. We answer their question for $k=3$ and sufficiently large $n$. We show that the unique minimum family is obtained by partitioning the ground set $[n]$ into two sets $A$ and $B$ with almost equal sizes and taking the family consisting of all the proper supersets of $A$ and of $B$.
\end{abstract}

\section{Introduction}

A central topic in extremal set theory is intersection properties of families of finite sets.  This topic began with classical results of Erd\H{o}s, Ko and Rado~\cite{EKR61} determining the maximum size of an intersecting family of subsets of an $n$-element set, for both arbitrary subsets and for subsets of a given size.

A family $\mF$ of sets is \emph{maximal} (or \emph{saturated}) with respect to some property if it satisfies the property, but no family properly containing $\mF$ satisfies the property.
The problem of finding the smallest object which is maximal with respect to some property has an extensive literature.  In the setting of graphs the topic was initiated by Erd\H{o}s, Hajnal and Moon~\cite{EHM}, who determined the minimum number of edges in a maximal $n$-vertex graph not containing a clique of size $k$ (see~\cite{survey} for an extensive survey on graph saturation). For set systems consisting of sets of a fixed size $r$, Erd\H{o}s and Lov\'asz~\cite{EL} suggested the problem of finding a maximal intersecting family of minimum size.
For this problem, the best known lower bound of $3r$ (when there is a projective plane of order $r-1$) was proved by Dow, Drake, F\"uredi and Larson~\cite{Dow}, and the best known upper bound of
$r^2/2 + O(r)$ is due to Boros, F\"uredi and Kahn~\cite{Boros}.  Kahn~\cite{kahn} conjectured that $O(r)$ is an upper bound.  Other properties for which saturation results have been obtained for set systems include the $k$-Sperner~\cite{Sp}, equal disjoint union free~\cite{DH}, VC-dimension at most $k$~\cite{VC} and $k$-matching free~\cite{Bucic} properties.

In this paper we study a natural saturation problem concerning an intersection property. To state the problem we need to introduce some notation. For a positive integer $n$, denote by $[n]$ the set $\{1,2,\ldots, n\}$. Given a set $S$, we write $2^S$ for the power set of $S$, which is 
the family of all subsets of $S$. %and $S^c$ for %the complement set of $S$, which is $[n]\sm S$.
For $\mF \subseteq 2^{[n]}$, let $\bar{\mF} \ce \{F^c : F \in \mF\}$ where $F^c$ denotes $[n]\sm F$. %Then, $\bar{\bar{\mF}} = \mF$ and $|\bar{\mF}| = |\mF|$.
We use the notation $\mF\Delta \mG$ for the symmetric difference $(\mF\sm \mG)\cup (\mG\sm \mF)$ of two set families $\mF$ and $\mG$.

We say that a family $\mF \subseteq 2^{[n]}$  is \emph{maximal $k$-wise intersecting} if 
the intersection of every collection of at most $k$ sets in $\mF$ is non-empty, and no set from $2^{[n]} \sm \mF$ can be added to $\mF$ while preserving this property. In the case $k=2$, it is well-known that every maximal $2$-wise intersecting family has the same size, namely $2^{n-1}$.  For $k \ge 3$, an old question of Erd\H{o}s
and Kleitman~\cite{EK74} from 1974 asks for the size of the smallest maximal $k$-wise intersecting family. We answer this question for the case $k=3$ and $n$ sufficiently large.

We now present our construction, which we call \emph{a balanced pair of linked cubes}.
A balanced pair of linked cubes is a set family of the form $\{A: S \subsetneq A \subseteq [n]\} \cup \{B: S^c \subsetneq B \subseteq [n]\}$ for some $S \subseteq [n]$ with $|S|=\lfloor n/2 \rfloor$.
It is not hard to check that a balanced pair of linked cubes is a maximal $3$-wise intersecting family of size $2^{\lfloor n/2 \rfloor} + 2^{\lceil n/2 \rceil} -3$.

Our main result is that for sufficiently large $n$, the smallest maximal $3$-wise intersecting families are balanced pairs of linked cubes.

\begin{theorem}\label{thm::main}
    If $\mF$ is a maximal $3$-wise intersecting family  of minimum size on ground set $[n]$, where $n$ is sufficiently large, then $\mF$ is a balanced pair of linked cubes.
\end{theorem}
\begin{remarkstar}
The even case of \cref{thm::main} was first proved by the third, fourth, sixth and seventh authors \cite{HLTT21}. In \cite{BCL22}, the first, second and fifth authors used the method of \cite{HLTT21} together with some new ideas to establish the odd case. The present paper provides a unified treatment of both cases.
\end{remarkstar}

%In order to complete the proof of \cref{thm::main}, 
A key ingredient in the proof of \cref{thm::main} is the following stability result, which shows that a maximal $3$-wise intersecting family of size at most $(1+o(1))(2^{\lfloor n/2\rfloor}+2^{\lceil n/2 \rceil})$ must be close in structure to a balanced pair of linked cubes.

\begin{theorem} \label{mlem}
	For every $\eps > 0$ there exists $\delta>0$ such that, if $\mF$ is a maximal $3$-wise intersecting family on $[n]$ of size $|\mF| \le (1+\delta) (2^{\lfloor n/2\rfloor}+2^{\lceil n/2 \rceil})$, where $n$ is sufficiently large, then there exists a balanced pair of linked cubes %$\mF_0 = \{A: A\subseteq S\}\cup \{B: B \subseteq S^c \}$, where $S\subseteq [n]$ with $|S|=\left\lfloor n/2 \right\rfloor$, 
	$\mF_0$ such that %$|\bar{\mF} \Delta \mF_{n,2}| \le \eps 2^{\lfloor n/2\rfloor}$.
$|\mF \tri \mF_0| \le \eps 2^{\lfloor n/2\rfloor}$.
\end{theorem}

%In other words, for a maximal $3$-wise intersecting family of size $(1+o(1))(2^{\lfloor n/2\rfloor}+2^{\lceil n/2 \rceil})$, all but a vanishing proportion of the members of the family are contained in a balanced pair of linked cubes.

\paragraph{Organization and notation:}
The rest of the paper is organized as follows. In \cref{sec:::stability} we prove \cref{mlem}, in \cref{sec::pro} we deduce \cref{thm::main} from \cref{mlem} using a perturbation argument, and in \cref{sec::largek} we discuss maximal $k$-wise intersecting families when $k\geq 4$.

We use standard asymptotic notation throughout. For functions $f = f(n)$ and $g = g(n)$ we write $f=o(g)$, $f\ll g$, or $g \gg f$ to mean $f/g \rightarrow 0$ as $n\rightarrow \infty$. All asymptotics are as
$n \rightarrow \infty$.

\section{Proof of \cref{mlem}}\label{sec:::stability}

\subsection{Proof Overview}

In the proof of \cref{mlem}, we find it more convenient to work with $\bar{\mF}$. 

\begin{definitionstar}
\emph{A balanced pair of cubes on $[n]$} is a set system of the form $\{A\colon A \subsetneq S\}\cup\{B\colon B\subsetneq S^c\}$, where $S\subseteq [n]$ with $|S|=\lfloor \frac{n}{2}\rfloor$. More generally, a \emph{balanced series of $k$ cubes on $[n]$}, denoted by $\mF_{n,k}$, 
is a set system of the form $\cup_{i=1}^k\{A\colon A \subsetneq S_i\}$, where $S_1,\ldots,S_k$ is a partition of $[n]$ with $\left\lfloor\frac{n}{k}\right\rfloor \le |S_i| \le  \left\lceil\frac{n}{k}\right\rceil$ for each $i\in[k]$.
\end{definitionstar}

We also need the following key concepts, whose relevance will be revealed in \cref{lem::toGen} below.

\begin{definitionstar}
	A set system $\mH\subseteq 2^{[n]}$ is a \emph{$(1-\eps)$-$k$-generator for $[n]$} if all but at most $\eps 2^n$ subsets of $[n]$ can be expressed as a union of at most $k$ disjoint sets of $\mH$.
\end{definitionstar}

\begin{definitionstar}
	Given a set system $\mH\subseteq 2^{[n]}$, denote by $\disj(\mH)$ the %set 
number of disjoint pairs in $\mH$, i.e., $\disj(\mH) =\left|\{\{A,B\}\colon A,B\in \mH,\, A\cap B=\emptyset\}\right|.$
\end{definitionstar}

The following lemma is central to the proof of \cref{mlem}.

\begin{lemma} \label{lem::toGen}
	If $\mF$ is a maximal $3$-wise intersecting family on $[n]$, then
	\begin{equation*}
		\disj(\bar{\mF} )\geq 2^n-|\mF|.
	\end{equation*} Moreover, if $|\mF| \le \eps 2^n$, then $\bar{\mF}$ is a $(1-\eps)$-$2$-generator for $[n]$.
\end{lemma}
\begin{proof}
	Let $\mF^c \ce 2^{[n]}\sm \mF$, then $|\mF^c|=2^n-|\mF|$. For every $A\in \mF^c$, there exist $B,C\in\mF$ such that $B\cap C\subseteq A^c$, since $\mF$ is a maximal $3$-wise intersecting family. Notice that $\mF$ is %closed upwards
upward closed, i.e.,~if $F\in \mF$ and $F\subseteq F'\subseteq [n]$, then $F'\in \mF$. 
Thus, we may select these sets $B,C\in \mF$ so that $A^c=B\cap C$ and $B\cup C=[n]$, which means $B^c\cap C^c=\emptyset$.
%Thus, we may choose $B\cap C=A^c$ with $B\cup C=[n]$, which is equivalent to $A=B^c\cup C^c$ with $B^c\cap C^c=\emptyset$. 
This (not necessarily unique) choice of $B$ and $C$ defines an injective map from $\mF^c$ into disjoint pairs of sets in $\bar{\mF}$ (where $\emptyset$ is mapped to $\{\emptyset, \emptyset\}$). Therefore, 
$$
	\disj(\bar{\mF} )\geq |\mF^c|=2^n-|\mF|.$$
	%In particular, every set in $\mF^c$ can be expressed as a union of at most two disjoint sets of $\bar{\mF}$. If $|\mF| \le \eps 2^n$, then $$|\mF^c|=2^n-|\mF| \geq 2^n- \eps 2^n = (1-\eps)2^n,$$ which implies that $\bar{\mF}$ is a $(1-\eps)$-$2$-generator for $[n]$.
	Now assume that $|\mF| \le \eps 2^n$. Since $|\mF^c|=2^n-|\mF| \geq (1-\eps)2^n$ and every set in $\mF^c$ can be expressed as a union of at most two disjoint sets of $\bar{\mF}$, $\bar{\mF}$ is a $(1-\eps)$-$2$-generator for $[n]$.
\end{proof}

The following result together with \cref{lem::toGen} %implies immediately 
immediately implies \cref{mlem}. We will prove it in \cref{subsec:ES}. 
%, by modifying ideas of Ellis and Sudakov~\cite{EllisSudakov2011}.\\

\begin{theorem}\label{thm::modify}
	For every $\eps'>0$, there exists $\delta=\delta(\eps')>0$ %and $\eta = \eta(\eps')>0$
such that for sufficiently large $n$, if $\mH\subseteq 2^{[n]}$ is a $(1-\delta)$-$2$-generator for $[n]$ with %$|\mH|\leq (1+\eta)(2^{\lfloor n/2\rfloor}+2^{\lceil n/2 \rceil})$
 $|\mH|\leq (1+\delta)(2^{\lfloor n/2\rfloor}+2^{\lceil n/2 \rceil})$, then there exists a balanced pair of cubes %$\mF_0 = \{A: A\subseteq S\}\cup \{B: B \subseteq S^c \}$, where $S\subseteq [n]$ with $|S|=\left\lfloor n/2 \right\rfloor$, 
	$\mF_{n,2}$ such that $|\mH\tri \mF_{n,2}|\leq \eps' 2^{\lfloor n/2\rfloor}$.
\end{theorem}

\begin{remarkstar}
\cref{thm::modify} for $n$ even was first proved by Ellis and Sudakov \cite[Proposition 9]{{EllisSudakov2011}}. Our proof of \cref{thm::modify} uses their method.    
\end{remarkstar}

\begin{proof}[Proof of \cref{mlem} assuming \cref{thm::modify}]
	 Let $\eps>0$ and $n$ be a sufficiently large integer. Let $\delta=\delta(\eps)>0$ be obtained from \cref{thm::modify}. %Let $\mF$ be a maximal $3$-wise intersecting family on $[n]$ of size $|\mF| \le (1+\eta) (2^{\lfloor n/2\rfloor}+2^{\lceil n/2 \rceil})$, then there exists some $\eps_1< \delta$ such that $|\mF|\leq \eps_1 2^n$, when $n$ is sufficiently large. 
	 Let $\mF$ be a maximal $3$-wise intersecting family on $[n]$ of size $|\mF| \le (1+\delta) (2^{\lfloor n/2\rfloor}+2^{\lceil n/2 \rceil})$, then $|\mF|\leq \delta 2^n$, when $n$ is sufficiently large. %Let $\mG\ce \bar{\mF}$, then $|\mG|=|\mF|\le (1+\eta) (2^{\lfloor n/2\rfloor}+2^{\lceil n/2 \rceil})$ and $\mG$ is a $(1-\eps_1)$-$2$-generator for $[n]$ by \cref{lem::toGen}. Thus, $\mG$ is a $(1-\delta)$-$2$-generator for $[n]$. 
    Since $|\bar{\mF}|=|\mF|\le \delta 2^n$, $\bar{\mF}$ is a $(1-\delta)$-$2$-generator for $[n]$ by \cref{lem::toGen}.
	 %By \cref{thm::modify}, there exists a balanced pair of cubes $\mF_0 = \{A: A\subseteq S\}\cup \{B: B \subseteq S^c \}$, where $S\subseteq [n]$ with $|S|=\left\lfloor n/2 \right\rfloor$, such that $|\bar{\mF}\tri\mF_0|=|\mG\Delta \mF_0|\leq \eps 2^{\lfloor n/2\rfloor}$.
	 Applying \cref{thm::modify} to $\mH:=\bar{\mF}$, we conclude that there exists a balanced pair of cubes $\mF_{n,2}$ with $|\bar{\mF}\tri\mF_{n,2}|\leq \eps 2^{\lfloor n/2\rfloor}$. Taking the complements of the sets in $\bar{\mF}$ and $\mF_{n,2}$ yields the desired result.
\end{proof}

\cref{lem::toGen} also leads to another classic question posed by Erd\H{o}s \cite{Erdos81}: How many disjoint pairs of sets can there be in a set system of given size? A lower bound comes from $\mF_{n,k}$, a balanced series of $k$ cubes, i.e. $\mF_{n,k}=\cup_{i=1}^k\{A: A \subsetneq S_i\}$, where $S_1,\ldots,S_k$ is a partition of $[n]$ with $\left\lfloor\frac{n}{k}\right\rfloor \le |S_i| \le  \left\lceil\frac{n}{k}\right\rceil$ for each $i\in[k]$. 
With the additional assumption that~$k$ divides~$n$, it is easy to see that $|\mF_{n,k}|= k2^{n/k}-2k+1$ and $\disj(\mF_{n,k}) > (1-1/k) \binom{|\mF_{n,k}|}{2}$.
Solving a conjecture of Daykin and Erd\H{o}s \cite{Guy83}, Alon and Frankl \cite{AF85} proved that %a balanced series of $k$ cubes 
$\mF_{n,k}$ has asymptotically the maximum number of disjoint pairs, given its size.

\begin{theorem}[Alon--Frankl]\label{thm:Alon-Frankl}
	For every positive integer $k$, there exists $\beta=\beta(k)>0$ such that, if $\mH$ is a set system on $[n]$ of size $m\ce |{\mH}| \ge 2^{(1/(k+1)+\eps)n}$, where $\eps>0$, then
	\[
	\disj(\mH) \le \left(1-\frac{1}{k}\right)\binom{m}{2}+O\left(m^{2-\beta\eps^2}\right).
	\]
\end{theorem}

We strengthen \cref{thm:Alon-Frankl} by proving a stability result, showing families with close to $\left(1-\frac{1}{k}\right)\binom{m}{2}$ disjoint pairs must %be close in structure to 
resemble a balanced series of $k$ cubes.

\begin{theorem}\label{thm:disjoint pairs stability}
	For every $\eps>0$ and integer $k\ge 2$, there exists an $\eta>0$ such that for every sufficiently large $n$, if a set system $\mH$ on $[n]$ of size $m\ce|\mH| \ge (1-\eta)k2^{n/k}$ has at least $(1-\frac{1}{k}-\eta)\frac{m^2}{2}$ disjoint pairs, then $k$ divides $n$, and %all but at most $\eps m$ sets in $\mF$ are contained inside a balanced series of $k$ cubes.
	there exists a balanced series of $k$ cubes $\mF_{n,k}$ with $|\mH\tri\mF_{n,k}| \le \eps 2^{\lfloor n/k\rfloor}$.
\end{theorem}

The proof of \cref{thm:disjoint pairs stability}, given in \cref{sec:1.6}, borrows some ideas of Alon and Frankl \cite{AF85}, and Alon, Das, Glebov and Sudakov~\cite{ADGS15}. 
The even case of  \cref{thm::modify} is proved using \cref{thm:disjoint pairs stability}.

\begin{proof}[Proof of \cref{thm::modify} for even $n$ assuming \cref{thm:disjoint pairs stability}]
 Let $\eta$ be as given by Theorem~\ref{thm:disjoint pairs stability} for $k=2$ and $\eps>0$, and let $\delta=\eta^2$.
	Let $\mH' \supseteq \mH$ such that $|\mH'| = \lfloor (1+\delta)2^{(n+2)/2} \rfloor$, and every element of $\mH' \sm \mH$ has size greater than $n/2$.
	This is possible because $n$ is large and the number of sets of size greater than $n/2$ is nearly $2^{n-1}$. 
 Since $\mH'\supseteq \mH$ and $\mH$ is a $(1-\delta)$-$2$-generator for $[n]$, we have
	\[
|\mH'|+\disj(\mH') \geq |\mH|+\disj(\mH) \geq (1-\delta)2^n.
 \]
 Hence $\disj(\mH') \ge (1-\delta)2^n-(1+\delta)2^{(n+2)/2} \ge (1-2\delta)2^n$ for $n$ sufficiently large.
	On the other hand,
	\[
 \left(\frac{1}{2} - \eta\right)\frac{|\mH'|^2}{2} \leq \left(\frac{1}{2} - \eta\right)(1+\eta)^2 2^{n+1} = (1-3\eta^2-2\eta^3) 2^n \le (1-3\delta)2^n.
 \]
	It follows that $\disj(\mH') \geq (1/2 - \eta)|\mH'|^2/2$, and the hypotheses of Theorem \ref{thm:disjoint pairs stability} are satisfied for $\mH'$.
	Theorem~\ref{thm:disjoint pairs stability} gives a balanced pair of cubes $\mF_{n,2}$ with $|\mH'\tri \mF_{n,2}|\leq \eps 2^{n/2}$.
	By the second property of $\mH'$, $\mH'\sm \mH$ is disjoint from $\mF_{n,2}$.
	Hence, $|\mH\tri \mF_{n,2}|\leq |\mH'\tri \mF_{n,2}|\leq \eps 2^{n/2}$. 
\end{proof}

Before starting the proofs of Theorems~\ref{thm::modify} and~\ref{thm:disjoint pairs stability}, we introduce the so-called disjointness graph, which is the main object that we analyze in Sections~\ref{sec:1.6} and~\ref{subsec:ES}.

\begin{definitionstar}
	For a set system $\mH$, the \emph{disjointness graph} $G_\mH$ is the graph with vertex set $\mH$ and edge set $\{\{A,B\}\subseteq \mH\colon A\cap B =\es\}$. For two (not necessarily disjoint) set families $\mH_1,\mH_2$, the \emph{disjointness bipartite graph} $G_{\mH_1, \mH_2}$ is the bipartite graph with classes $(\mH_1,\mH_2)$, where there is an edge between $A \in \mH_1$ and $B \in \mH_2$ if and only if $A \cap B = \es$.
\end{definitionstar}

\subsection{Proof of \cref{thm:disjoint pairs stability}}\label{sec:1.6}

The heart of the proof of \cref{thm:disjoint pairs stability} is the following lemma, which shows that the disjointness graph of a large family has only few cliques of size $k+1$. This lemma essentially appears in \cite{AF85}. For completeness, we include its proof here.

\begin{lemma}\label{lem:cliques count}
	For every $\eps>0$, $\gamma>0$ and integer $k\ge 2$, and sufficiently large integer $n$, if $\mH$ is a set system on $[n]$ of size %$m\ce |\mF| \ge 2^{n/k}$
$m\ce |\mH| \ge 2^{(1/(k+1)+\eps)n}$, then $G_{\mH}$ contains at most $\gamma \binom{m}{k+1}$ copies of $K_{k+1}$.
\end{lemma}

We denote by $K_r(t)$ the complete $r$-partite graph with parts of size $t$. The following proposition is standard and follows from a result of Erd\H{o}s~\cite{Erdos64} by a simple averaging argument (see for instance \cite[Proposition~3.2]{ADGS15}).

\begin{proposition}\label{prop:blow up}
	For integers $r\ge 2$, $t \ge 1$ and any real $\gamma>0$, there exists $\delta_{\ref{prop:blow up}}=\delta_{\ref{prop:blow up}}(r,t,\gamma)>0$, such that if $m$ is sufficiently large and $G$ is a graph on $m$ vertices with at least $\gamma \binom{m}{r}$ copies of $K_r$, then $G$ contains at least $\delta_{\ref{prop:blow up}} \binom{m}{rt}$ copies of $K_r(t)$.
\end{proposition}

We will use a probabilistic argument to derive \cref{lem:cliques count} from \cref{prop:blow up}.

\begin{proof}[Proof of \cref{lem:cliques count}]
	%Let $t=2k(k+1)$. 
	Let $t=\lceil\frac{2}{\eps}\rceil+k$. Select $t$ sets $A_1,\ldots,A_t \in {\mH}$ independently uniformly at random, with repetitions allowed. The probability that $|A_1\cup A_2\cup \cdots \cup A_t| \le \frac{n}{k+1}$ is bounded above by
	\begin{equation*}
		\sum_{\substack{S\subseteq [n] \\ |S|=\lfloor n/(k+1)\rfloor} } \Pr\big[ A_i\subseteq S, i=1,\ldots,t\big]
		\le 2^n \left(2^{n/(k+1)}/m\right)^t \le m^{-k},
	\end{equation*}
	since
	$m\ge 2^{(1/(k+1)+\eps)n}$, where $2^n$ estimates the number of choices of $S$, and $S$ has $2^{\lfloor n/(k+1)\rfloor}$ many subsets, each could be chosen as $A_i$, if they are in $\mH$.
	
	Let $\mathbf{A}=(A_1,A_2,\ldots,A_{(k+1)t})$ be a random sample (chosen independently uniformly at random, allowing repetition) of $(k+1)t$ vertices of $G_{\mH}$. It follows from the discussion above that the probability that $|A_{(i-1)t+1}\cup A_{(i-1)t+2}\cup \cdots \cup A_{it}| \le \frac{n}{k+1}$ for some $i\in [k+1]$ is at most $(k+1)m^{-k}$. On the other hand, if our random sample $\mathbf{A}$ gives a copy of $K_{(k+1)}(t)$ with $\{A_{(i-1)t+1},\ldots, A_{it}\}$ being the $i$-th vertex class for every $i\in [k+1]$, then $A_{1}\cup \cdots \cup A_{t}, \ldots, A_{kt+1}\cup \cdots \cup A_{(k+1)t}$ are $k+1$ disjoint subsets of $[n]$. Hence, in this case, we must have $|A_{(i-1)t+1}\cup A_{(i-1)t+2}\cup \cdots \cup A_{it}| \le \frac{n}{k+1}$ for some $i\in [k+1]$. Moreover, the probability that our random sample gives such a copy of $K_{k+1}(t)$ is precisely
	\[
	m^{-(k+1)t}\cdot (k+1)! (t!)^{k+1} \cdot \# \ \text{copies of } K_{k+1}(t) \text{ in $G_\mathcal{\mH}$}.
	\]
	Therefore, the number of copies of $K_{k+1}(t)$ in $G_{\mH}$ is at most
	\[
	\frac{m^{(k+1)t-k}}{k! (t!)^{k+1}} \le \delta_{\ref{prop:blow up}}(k+1,t,\gamma)\binom{m}{(k+1)t}
	\]
	for $m$ sufficiently large. Hence, by \cref{prop:blow up}, $G_{\mH}$ has at most $\gamma \binom{m}{k+1}$ copies of $K_{k+1}$. This completes the proof of \cref{lem:cliques count}.
\end{proof}

We will make use of the following theorem of Balogh, Bushaw, Collares, Liu, Morris and Sha\-rif\-zadeh \cite[Theorem 1.2]{BBCLMS17}.

\begin{theorem}\label{thm:removal}
	For every $m,k,t \in \bN$, the following holds.  Suppose $G$ is graph on $m$ vertices with at most $\frac{m^{k-1}}{e^{2k}\cdot k!}\left(e(G)+t-\left(1-\frac{1}{k}\right)\frac{m^2}{2}\right)$ copies of $K_{k+1}$.
	Then there is a partition of the vertex set of $G$ as $V(G)=V_1\cup V_2\cup
	\cdots \cup V_k$ with $\sum_{i=1}^k e(V_i) \le t$.
\end{theorem}

We also use the following well-known estimate on the size of a set system in terms of the binary entropy function.
\begin{lemma}\label{lem:entropy}
	Let $\mF$ be a set system on $[n]$. Denote by $p_i$ the fraction of sets in $\mF$ that contain $i$. Then
	\[
	|\mF| \le 2^{\sum_{i=1}^n h(p_i)},
	\]
	where $h(p)=-p\log_2 p-(1-p)\log_2 (1-p)$ is the binary entropy function.
\end{lemma}

We are now ready to prove \cref{thm:disjoint pairs stability}. The proof closely follows the approach from~\cite{ADGS15}.

\begin{proof}[Proof of \cref{thm:disjoint pairs stability}]
	
	Let $\mH$ be a family of subsets of $[n]$ such that $m\ce |{\mH}| \ge (1 - o(1)) k 2^{n/k}$ and $\disj(\mH) \ge \left( 1-1/k-o(1) \right) \frac{m^2}{2}$.
	By \cref{lem:cliques count}, $G_{\mH}$ contains at most $o(m^{k+1})$ copies of $K_{k+1}$. Thus, applying \cref{thm:removal} with $t=o(m^2)$, we conclude that $G_{\mH}$ has a $k$-partite subgraph with $(1 - 1/k - o(1))\frac{m^2}{2}$ edges. In order to contain this many edges, it is clear that the vertex classes must all have size $(1 - o(1)) m/k > (1 - o(1)) 2^{n/k}$.
	Furthermore, again by the edge density, all but at most $o(m)$ vertices must have at least $(1-1/k-o(1))m$ neighbors.
	Thus, we may take a $k$-partite subgraph $G$ of $G_{\mH}$ consisting of parts of size at least $(1 - o(1)) 2^{n/k}$ such that each vertex has
	at least $(1-1/k-o(1)) k 2^{n/k}$ neighbors, and thus at least $(1-o(1)) 2^{n/k}$ neighbors in every other class.

	Let $\eps>0$ be a sufficiently small constant with respect to $k$. Let $\mH_1, \mH_2,\dots,\mH_k$ be the sets from $\mH$ corresponding to the $k$ color classes of $G$, respectively. By the above discussion,
	\begin{equation}\label{lower}
		|\mH_i|\ge (1-\eps)2^{n/k} \quad \text{for all} \enskip i\in [k].
	\end{equation}
	We will construct a balanced series of $k$ cubes containing $\bigcup_{i=1}^k \mH_i$.  To do this we begin with a partition of the ground set $[n]$ into $2k+1$ disjoint sets: $A_1,B_1,\ldots,A_k,B_k,R$.  For $i \in [k]$, let $A_i$ and $B_i$ be the set of those elements of $[n] \sm \bigcup_{j=1}^{i-1} (A_j \cup B_j)$ occurring in more than $\eps |\mH_i|$ members of $\mH_i$ and between $1$ and $\eps |\mH_i|$  members of $\mH_i$, respectively. Let $R$ be the set of remaining elements of $[n]$.

	By the definition of the partition, if for some $i \in [k]$ we have $x \in B_i$, then there is a set $F \in \mH_i$ such that $x \in F$. It follows from the density properties of $G$ discussed above that for $j \neq i$, $x$ is not contained in at least $(1 - \eps) |{\mH_j}|$ members of $\mH_j$.
	If $x \in A_i$, then for $j\neq i$, $x$ does not appear in any member of $\mH_{j}$, since each member of $\mH_{j}$ is disjoint from at least $(1 - \eps) |\mH_i|$ members of $\mH_i$.
	
	It follows that  no set $F\in \mH$ has vertices from $\bigcup_{j=1}^{i-1} A_j$.  Furthermore, each element $x \in \bigcup_{j=1}^{i} B_j$ is contained in at most $\eps |\mH_i|$ sets in $\mH_i$.
	Using \cref{lem:entropy} we have
	\begin{equation*}
		|\mH_i|\leq 2^{|A_i|+h(\eps)(\sum_{j \le i}|B_j|)}.
	\end{equation*}
	We thus get
	\begin{align*}
		2^{n-\frac{1}{2}} < (1-\eps)^k 2^{n}&\leq \prod_{i=1}^k|\mH_i| \le 2^{\sum_{i\in [k]}\left(|A_i| + h(\eps)(k-i+1)|B_i|\right)}\\
		&= 2^{n - |R|-\sum_{i\in [k]} \left(1-h(\eps)(k-i+1)\right) |{B_i}|} x\le 2^{n -|R|- \frac{3}{4}\sum_{i\in [k]}|B_i|},
	\end{align*}
	where in the first equality we used the fact that $n=|R|+\sum_{i\in [k]}(|A_i|+|B_i|)$, and the last inequality holds since $1-h(\eps)k \ge 3/4$ for $\eps$ sufficiently small with respect to $k$. This implies that $R=B_1=\cdots=B_k=\emptyset$. As $\mH_i \subseteq 2^{A_i\cup B_i}$, by the definition of $A_i$ and $B_i$, we have that $|\mH_i| \le 2^{|A_i|+|B_i|}=2^{|A_i|}$. Together with the lower bound $|\mH_i|\ge (1-\eps)2^{n/k}$ from \eqref{lower}, we must have $|{A_i} > (n-1)/k$, and so $\card{A_i} \ge n/k$.  Since $\sum_{i \in [k]} \card{A_i} = n$, we in fact have $|A_i|=n/k$ for all $i\in [k]$.  It follows that the family $\bigcup_{i=1}^{k} \mH_i$ is contained in the balanced series of $k$ cubes $A_1,\ldots,A_k$.  Since $\mH$ contains at most $\eps |\mH|$ members not in  $\bigcup_{i=1}^{k} \mH_i$, the proof is complete.
\end{proof}

\subsection{Proof of \cref{thm::modify}}\label{subsec:ES}

Before stating the proof of \cref{thm::modify}, we need some preparation.

\begin{definitionstar}
	Given set systems $\mH_1, \mH_2$ and a bipartite subgraph $B\subseteq G_{\mH_1}$ with bipartition $(\mX,\mY)$, we say $B$ (and sometimes we say $E(B)$) \emph{generates} $\mH_2$ if every $H\in \mH_2$ can be expressed as a disjoint union of some $X\in\mX$ and $Y\in\mY$, i.e., every $H\in \mH_2$ corresponds to an edge of $B$.
\end{definitionstar}

\begin{definitionstar}
	For a set system $\mH$ and $i\in [n]$, let $\mH_i^- \ce \{H\in \mH: i\notin H\}$ be the subfamily of $\mH$, consisting of sets not containing $i$ and $\mH_i^+ \ce \{H\sm\{i\}: i\in H\in \mH\}$. Note that $|\mH_i^+|+|\mH_i^-|=|\mH|$.
\end{definitionstar}

We will use the following result by Ellis and Sudakov \cite{EllisSudakov2011}.

\begin{lemma}[{\cite[Proposition 18]{EllisSudakov2011}}]\label{lem::bipartite}
	Let $c>0$. Then, %there exists $b=b(c)>0$ such that 
 for every set family $\mH\subseteq 2^{[n]}$ with $|\mH|\geq c2^{n/2}$, the disjointness graph $G_{\mH}$ can be made bipartite by deleting at most %$\frac{(\log\log n)^b}{\log n}|\mH|^2$
 $o(|\mH|^2)$ edges.
\end{lemma}

%%%%%%%%%%%%%%%%%%%%%
%%%%%%%%%%%%%%%%%%%%%
\begin{comment}

\begin{lemma}[{\cite[Proposition 18]{EllisSudakov2011}}]\label{lem::bipartite}
	Let $c>0$. Then, there exists $b=b(c)>0$ such that for every $\mH\subseteq 2^{[n]}$ with $|\mH|\geq c2^{n/2}$, the disjointness graph $G_{\mH}$ can be made bipartite by deleting at most $\frac{(\log\log n)^b}{\log n}|\mH|^2$ edges.
\end{lemma}

In~\cite{EllisSudakov2011}, Ellis and Sudakov proved \cref{thm::modify} for even $n$, which will be needed for the proof of \cref{thm::modify}. We state it below.

\begin{theorem}[{\cite[Proposition 9]{EllisSudakov2011}}]\label{thm::stability n even}
	For every $k\in\mathbb{N}$ and $\eps'>0$, there exist $n_0=n_0(k,\eps')$, $\delta = \delta(k, \eps')>0$ and $\eta=\eta(k,\eps')>0$ such that the following holds. If $n\geq n_0$ is a multiple of $k$ and $\mG\subseteq 2^{[n]}$ is a $(1-\delta)$-$k$-generator for $[n]$ with %$|\mG|\leq (1+\eta) |\mF_{n,k}|$, 
	\textcolor{red}{$|\mG|\leq (1+\eta) (k2^{n/k}-k+1)$}, then there exists a balanced series of $k$ cubes %$\mF_1$ such that $|\mG\tri \mF_1|\leq \eps' |\mF_{n,k}|$.
	$\mF_{n,k}$ such that $|\mG\tri \mF_{n,k}|\leq \eps' |\mF_{n,k}|$.
\end{theorem}

\end{comment}

%%%%%%%%%%%%%%%%%%%%%%%%%%
%%%%%%%%%%%%%%%%%%%%%%%%%%%%%%

We are now ready to prove \cref{thm::modify}.

\begin{proof}[Proof of \cref{thm::modify}.]
	%By \cref{thm::stability n even}, we may assume that $n=2\ell+1$ is a sufficiently large odd integer. 
	We have already proved the even case of \cref{thm::modify}, hence  from now on we assume that $n=2\ell+1$ is a sufficiently large odd integer.
	%Let $\delta \gg \etan >0$ be sufficiently small.
 Let $\eps'\gg\eps \gg \delta>0$ be sufficiently small. Suppose that $\mH \subseteq 2^{[n]}$ is a $(1-\delta)$-$2$-generator for $[n]$ with $|\mH|\leq (1+\etan) (2^{\lfloor n/2\rfloor}+2^{\lceil n/2 \rceil})=3(1+\etan)2^\ell$. Then, the number of ways to choose at most two disjoint sets (whose unions are all different) from $\mH$ is at least $(1-\delta)2^n$, by the definition of a $(1-\delta)$-$2$-generator. Hence, $|\mH|^2\geq (1-\delta)2^n$, implying %which implies that 
	\[|\mH|\geq \sqrt{1-\delta}\cdot 2^{n/2}.\]
	Moreover, let %$G_0 \ce G_{\mH}$ 
	$G_0$ be the disjointness graph of $\mH$, then \[
	1+|\mH|+e(G_0)\geq (1-\delta)2^n,\] which implies that \[
	e(G_0)\geq (1-\delta)2^{2\ell+1}-3(1+\etan)2^\ell-1.\] We conclude that $G_0$ has edge-density
	$$\frac{e(G_0)}{\binom{|\mH|}{2}}\geq \frac{(1-\delta)2^{2\ell+1}-3(1+\etan) 2^\ell-1}{9(1+\etan)^2 2^{2\ell-1}} \geq \frac{4-5\delta}{9(1+\etan)^2}, $$ where the last inequa\-lity comes from $3(1+\etan)2^\ell+1\leq \delta2^{2\ell-1}$.
	Applying \cref{lem::bipartite} to $\mH$ with $c=\sqrt{1-\delta}$, we conclude that one can delete at most %we get that there exists a constant $b>0$ such that we can delete at most
	%$$
	%\frac{(\log\log n)^b}{\log n}|\mH|^2\leq \frac{(\log\log n)^b}{\log n}9(1+\etan)^22^{2\ell}=\frac{9(1+\etan)^2(\log\log n)^b}{2\log n}2^{2\ell+1}
	%$$
	$$
	o(|\mH|^2)=o(2^{2\ell+1})
	$$
	edges from $G_0$ and obtain a bipartite graph $G=(\mX,\mY)$ with $\mX\cup\mY=\mH$. Note that $G$ generates all but at most
	%$$
	%\delta2^{2\ell+1}+1+|\mH|+\frac{9(1+\etan)^2(\log\log n)^b}{2\log n}2^{2\ell+1}
	%$$
	$$
	\delta2^{2\ell+1}+1+|\mH|+o(2^{2\ell+1})
	$$
	subsets of $[n]$. %Since $\frac{9(1+\etan)^2(\log\log n)^b}{2\log n}=o(1)$ and $|\mH|\leq 3(1+\etan) 2^\ell=o(2^{2\ell+1})$, we may assume that $G$ generates all but at most $\eps2^n$ subsets of $[n]$, where $\eps=2\delta = o(1)$ satisfies
	%\begin{equation}\label{eG}
	%	e(G)\geq (1-\eps)2^{2\ell+1}.	
	%\end{equation}
    %In particular, $\mH$ is a $(1-\eps)$-$2$-generator for $[n]$.
    Since $|\mH|\leq 3(1+\etan) 2^\ell=o(2^{2\ell+1})$ and $\etan \ll \eps$, $G$ generates all but at most $\eps 2^{2\ell+1}$ subsets of $[n]$. %, where $\eps:=2\delta =o(1)$. 
    In particular,
    \begin{equation}\label{eG}
		e(G)\geq (1-\eps)2^{2\ell+1}.	
	\end{equation}

    Let $\alpha \ce |\mX|/2^\ell$ and $\beta \ce |\mY|/2^\ell$. Since $|\mX| + |\mY| = |\mH| \le 3 (1+\etan) 2^\ell$, we have
    \begin{equation} \label{equ::alpplusbet}
    	\alpha+\beta\leq 3(1+\etan).
    \end{equation}
    Therefore, we have
    \begin{equation} \label{equ::alpmulbet}
    	\alpha\beta\leq \frac{9}{4}(1+\etan)^2.
    \end{equation}
    Since $\alpha \beta 2^{2\ell}=|\mX||\mY|\geq e(G)\geq (1-\eps)2^{2\ell+1}$, we also have $\alpha\beta\geq 2-2\eps$. Combining with \eqref{equ::alpplusbet}, we get
    \begin{equation}\label{eqn::bound alphabeta}
    	1-3\eps<\alpha,\beta<2+3\eps,
    \end{equation}
    where we use %$\eps= 2\delta \gg \etan$
    $\eps \gg \etan$.

    Let
    \begin{equation}\label{1/3def}
    	\mX_{1/3}\ce\{i\in [n]: |\mX_i^+|\geq |\mX|/3\}\quad \text{\ \ \ and\ \ \ \ \ \ }
    	\mY_{1/3}\ce\{i\in [n]: |\mY_i^+|\geq |\mY|/3\}.
    \end{equation}
    Letting $A\ce\mX_{1/3}$ and $B\ce\mY_{1/3}$, we will show that $A, B$ forms an equipartition of $[n]$, and that $\mX$ and $\mY$ are not too far from $2^A$ and $2^B$, respectively, thereby proving \cref{thm::modify}. %We characterize the structure of $A\ce\mX(1/3)$ and $B\ce\mY(1/3)$ step by step via the following lemma and a series of claims.
    We will achieve this goal step by step via the following lemma and a series of claims.

    \begin{lemma}\label{lem::size}
    	It cannot happen that all the following equations hold at the same time:
    	\begin{equation}\label{eqn::size}
    		|\mX|, |\mY|=(3/2-o(1))2^\ell,\quad
    		|\mX_n^+|, |\mY_n^+|=(1-o(1))2^{\ell-1} \quad  \textrm{and} \quad
    		|\mX_n^-|, |\mY_n^-|=(1-o(1))2^\ell.
    	\end{equation}
    \end{lemma}

\begin{proof}%of Lemma~\ref{lem::size}
	Suppose for a contradiction that~\eqref{eqn::size} holds. Since all but at most $\eps2^n=2\eps \cdot 2^{n-1}$ sets in $[n]$ can be expressed as a union of at most two disjoint sets in $\mH$, we have $\mH_n^-=\mX_n^- \cup \mY_n^-$ is a $(1-2\eps)$-$2$-generator for $[n-1]$. By the assumption that~\eqref{eqn::size} holds, $|\mH_n^-|=(2-o(1))2^\ell=(2-o(1))2^{(n-1)/2}$. As $|\mF_{n-1,2}|=2\cdot 2^{(n-1)/2}-1$, we can apply %\cref{thm::stability n even} 
	the even case of \cref{thm::modify} to $\mH_n^-$ and conclude that there exists an equipartition $S_1\cup S_2=[n-1]$ such that each of $\mX_n^-, \mY_n^-$ contains at least $(1-o(1))2^\ell$ sets in $2^{S_1}, 2^{S_2}$. Define $\mU\ce\{F\in \mX: F\cap S_2=\emptyset\}, \mV\ce\{F\in \mY: F\cap S_1=\emptyset\}$.
	Note that $\mU_n^- = \mX_n^- \cap 2^{S_1}$ and $\mV_n^- = \mY_n^- \cap 2^{S_2}$.
	We have $|\mU_n^-|, |\mV_n^-|=(1-o(1))2^\ell$, implying
	\begin{equation}\label{eqn::size minus}
		|\mX_n^-\sm \mU_n^-|, |\mY_n^-\sm \mV_n^-|=o(2^\ell).
	\end{equation}
	
	\noindent Now we prove that \begin{equation*}
		|\mX_n^+\sm \mU_n^+|, |\mY_n^+\sm \mV_n^+|=o(2^\ell).
	\end{equation*} In fact, for every $X\in\mX_n^+\sm \mU_n^+$, we have $X\cap S_2\neq \emptyset$ by the definition of $\mU$, thus $X\cup\{n\}$ is disjoint from at most $2^{\ell-1}$ subsets of $S_2$. Since $|\mY_n^-\sm 2^{S_2}|=|\mY_n^-\sm \mV_n^-|=o(2^\ell)$, the set $X\cup\{n\}$ is disjoint from at most $2^{\ell-1}+o(2^\ell)$ sets in $\mY$. Similarly, for every $Y\in\mY_n^+\sm \mV_n^+$, the set $Y\cup\{n\}$ is disjoint from at most $2^{\ell-1}+o(2^\ell)$ sets in $\mX$. Let $e_n$ be the number of edges $XY\in E(G)$ such that $n\in X\cup Y$. Since $G$ generates all but at most $\eps2^n$ subsets of $[n]$, at least $2^{2\ell}-\eps2^n=(1-2\eps)2^{2\ell}$ sets containing $n$ correspond to edges of $G$, which implies that \begin{equation*}
		e_n\geq (1-2\eps)2^{2\ell}=\left(1-o(1)\right)2^{2\ell}.
	\end{equation*} Let $\phi\ce |\mU_n^+|/|\mX_n^+|$ and $\theta\ce |\mV_n^+|/|\mY_n^+|$, then $\phi,\theta\in[0,1]$. Combining with~\eqref{eqn::size}, we have \begin{equation*}
		\begin{aligned}
			e_n &\leq |\mU_n^+||\mY_n^-|+|\mX_n^+\sm\mU_n^+|\left(2^{\ell-1}+o(2^\ell)\right)+|\mV_n^+||\mX_n^-|+|\mY_n^+\sm\mV_n^+|\left(2^{\ell-1}+o(2^\ell)\right)\\
			&=\phi(1-o(1))2^{2\ell-1}+(1-\phi)(1-o(1))2^{2\ell-2}+\theta(1-o(1))2^{2\ell-1}+(1-\theta)(1-o(1))2^{2\ell-2}\\
			&=\left(2+\phi+\theta-o(1)\right)2^{2\ell-2}.
		\end{aligned}
	\end{equation*} Hence, 
    \begin{equation*}
		\left(1-o(1)\right)2^{2\ell}\leq e_n\leq \left(2+\phi+\theta-o(1)\right)2^{2\ell-2},
    \end{equation*} 
    which implies that $\phi,\theta=1-o(1)$. Therefore, $|\mX_n^+\sm \mU_n^+|=o(|\mX_n^+|)=o(2^{\ell})$, and $|\mY_n^+\sm \mV_n^+|=o(|\mY_n^+|)=o(2^{\ell}),$ as desired.
	
	Now we have \begin{equation*}
		|\mX\sm \mU|+ |\mY\sm \mV|=|\mX_n^-\sm \mU_n^-|+|\mX_n^+\sm \mU_n^+|+|\mY_n^-\sm \mV_n^-|+|\mY_n^+\sm \mV_n^+|=o(2^\ell). \end{equation*}
	For a set $F\in\mH$, notice that $F\in (\mX\sm \mU)\cup (\mY\sm \mV)$ if and only if $F\in\mX$ satisfies $F\cap S_2\neq \emptyset$ or $F\in\mY$ satisfies $F\cap S_1\neq \emptyset$. For $S\subseteq S_1$, if $S\cup \{n\}\in \mX$, then $S\in\mX_n^+$ by definition. Hence, there are at least $2^{|S_1|}-|\mX_n^+|=2^\ell-(1-o(1))2^{\ell-1}=(1+o(1))2^{\ell-1}$ sets $S\subseteq S_1$ satisfying $S\cup \{n\}\notin \mX$. Similarly, there are at least $(1+o(1))2^{\ell-1}$ sets $S'\subseteq S_2$ satisfying $S'\cup \{n\}\notin \mY$. Thus, there are at least $(1+o(1))2^{2\ell-2}$ sets of the form $S\cup S'\cup\{n\}\subseteq [n]$ satisfying $S\subseteq S_1, S\cup \{n\}\notin \mX$ and $S'\subseteq S_2, S'\cup \{n\}\notin \mY$. 
	Denote the family of sets of this form by $\mS$. Since $G$ generates all but at most $\eps2^n$ subsets of $[n]$, at least $(1+o(1))2^{2\ell-2}-\eps2^{2\ell+1}=(1-o(1))2^{2\ell-2}$ sets in $\mS$ correspond to edges of $G$. If $F\in\mS$ can be expressed as a disjoint union of $X\in \mX$ and $Y \in \mY$, then either $X\cap S_2\neq\emptyset$ or $Y\cap S_1\neq \emptyset$, implying that either $X$ or $Y$ is in $(\mX\sm \mU)\cup (\mY\sm \mV)$. Hence, the number of choices for $F$ is at most $o(2^\ell)|\mH|\leq o(2^\ell)\cdot 3(1+\etan)2^\ell\ll(1-o(1))2^{2\ell-2}$, a contradiction.
\end{proof}

\begin{claim}\label{clm::union}
	$A\cup B=[n]$.
\end{claim}

\begin{proof}%of Claim~\ref{clm::union}
	Suppose for a contradiction that $A\cup B\neq [n]$. We may assume without loss of generality that $n\notin A\cup B$.
	%then $n$ is contained in less than $|\mX|/3$ sets of $\mX$ and less than $|\mY|/3$ sets of $\mY$.
	Let $x\ce|\mX_n^+|/|\mX|$ and $y\ce|\mY_n^+|/|\mY|$, then $x,y<1/3$ by the definitions of $A$ and $B$. Recalling that $e_n\geq (1-2\eps)2^{2\ell}$ is the number of disjoint pairs $X\in\mX, Y\in\mY$ such that $n\in X\cup Y$, we have \begin{equation}\label{eqn::derive size}
		\begin{aligned}
			(1-2\eps)2^{2\ell}\leq e_n&\leq |\mX_n^+||\mY_n^-|+|\mY_n^+||\mX_n^-|= |\mX_n^+|\left(|\mY|-|\mY_n^+|\right)+|\mY_n^+|\left(|\mX|-|\mX_n^+|\right)\\
			&=
			x\alpha2^\ell(\beta 2^\ell - y\beta 2^\ell) + y\beta2^\ell(\alpha2^\ell - x\alpha 2^\ell)
			=(x+y-2xy)\alpha\beta 2^{2\ell}.
		\end{aligned}
	\end{equation}
	
	Define the function $f(x,y)\ce x+y-2xy$. On $0\leq x,y\leq 1/3$, the function $f(x,y)$ attains its maximum value $4/9$, when $x=y=1/3$. Combining with~\eqref{eqn::derive size}, we have
	$$1-2\eps\leq f(x,y)\alpha\beta\leq \frac{4}{9}\alpha\beta.$$
	Recalling that $\alpha+\beta\leq 3(1+\etan)$ by \eqref{equ::alpplusbet}, we then get \begin{equation*}
		\frac{3}{2}-4\sqrt{\frac{\eps}{2}}\leq \alpha,\beta\leq \frac{3}{2}+4\sqrt{\frac{\eps}{2}}.
	\end{equation*}
	Additionally, $\alpha\beta\leq \frac{9}{4}(1+\etan)^2$ by \eqref{equ::alpmulbet}, implying that $1-2\eps\leq f(x,y)\alpha\beta\leq \frac{9}{4}(1+\etan)^2f(x,y)$, so
	\begin{equation*}
		\frac{1}{3}- 3\eps \leq x,y\leq \frac{1}{3}.
	\end{equation*}
	
	In summary, we have \begin{equation*}
		|\mX|, |\mY|=(3/2-o(1))2^\ell,\quad
		|\mX_n^+|, |\mY_n^+|=(1-o(1))2^{\ell-1} \quad  \textrm{and} \quad
		|\mX_n^-|, |\mY_n^-|=(1-o(1))2^\ell,
	\end{equation*} where $\eps=o(1)$.
	Therefore, \cref{clm::union} follows from \cref{lem::size}.
\end{proof}

\begin{claim}\label{clm::intersection}
	$A\cap B=\emptyset$.
\end{claim}

\begin{proof}%of Claim~\ref{clm::intersection}
	Suppose for a contradiction that $A\cap B\neq \emptyset$. We may assume without loss of generality that $n\in A\cap B$. Let $x\ce|\mX_n^+|/|\mX|$ and $y\ce|\mY_n^+|/|\mY|$, then $x,y\geq 1/3$ by the definitions of $A$ and $B$. Notice that \begin{equation*}
		(2-2\eps)2^{2\ell}\leq e(G)\leq |\mX||\mY|-|\mX_n^+||\mY_n^+|=(1-xy)\alpha\beta 2^{2\ell}.
	\end{equation*}
	Hence,
	$$
	2-2\eps\leq (1-xy)\alpha\beta\leq \frac{8}{9}\alpha\beta.
	$$
	Recalling that $\alpha + \beta \le 3(1+\etan)$ by \eqref{equ::alpplusbet}, we then get
	\begin{equation*}
		\frac{3}{2}-2\sqrt{\eps}\leq \alpha,\beta\leq \frac{3}{2}+2\sqrt{\eps}.
	\end{equation*}
	Additionally, $\alpha\beta\leq \frac{9}{4}(1+\etan)^2$ by \eqref{equ::alpmulbet}, implying that $2-2\eps\leq (1-xy)\alpha\beta\leq \frac{9}{4}(1+\etan)^2(1-xy)$, so
	\begin{equation*}
		\frac{1}{3}\leq x,y\leq \frac{1}{3}+3\eps.
	\end{equation*}
	
	In summary, we have \begin{equation*}
		|\mX|, |\mY|=(3/2-o(1))2^\ell,\quad
		|\mX_n^+|, |\mY_n^+|=(1-o(1))2^{\ell-1} \quad  \textrm{and} \quad
		|\mX_n^-|, |\mY_n^-|=(1-o(1))2^\ell,
	\end{equation*} where again $\eps=o(1)$. By \cref{lem::size}, we have completed the proof of \cref{clm::intersection}.
\end{proof}

By Claims~\ref{clm::union} and~\ref{clm::intersection}, we now know that $A\cup B$ is a partition of $[n]$. It remains to show that $A\cup B$ is in fact an equipartition of $[n]$ and $\mX, \mY$ are close to $2^A,2^B$, respectively. %, thus we prove \cref{thm::modify}. 
The following observation is simple but will be useful hereafter.

\begin{observation}\label{obs}
	If $F\in \mX\sm 2^A$, then $F$ has at most $2|\mY|/3$ neighbors in $\mY$. Similarly, if $F\in \mY\sm 2^B$, then $F$ has at most $2|\mX|/3$ neighbors in $\mX$.
\end{observation}

\begin{proof}
	By symmetry, it suffices to prove the first part. Suppose $F\in \mX\sm 2^A$, then there exists $i\in [n]$ such that $i\in F\cap B$. By the definition of $B=\mY_{1/3}$ (see~\eqref{1/3def}), there are at least $|\mY|/3$ sets in $\mY$ containing $i$, which therefore have non-empty intersection with $F$.
	By the definition of $G$, we conclude that $F$ has at most $|\mY|-|\mY|/3=2|\mY|/3$ neighbors in $\mY$, as desired.
\end{proof}

\begin{claim}\label{clm::size}
	%We have $|\mX\cap 2^A|\geq (2/3-3\eps)|\mX|=(2/3-o(1))|\mX|$ and $|\mY\cap 2^B|\geq (2/3-3\eps)|\mY|=(2/3-o(1))|\mY|$.
We have $|\mX\cap 2^A|\geq (2/3-3\eps)|\mX|$ and $|\mY\cap 2^B|\geq (2/3-3\eps)|\mY|$. Additionally, $[n]=A\cup B$ is an equipartition.
\end{claim}

\begin{proof}%of Claim~\ref{clm::size}
	Let $\theta \ce \frac{|\mX \cap 2^A|}{|\mX|}$ and $\phi \ce \frac{|\mY \cap 2^B|}{|\mY|}$. By \cref{obs} and %$\alpha\beta\leq \frac{9}{4}(1+\etan)^2$
	\eqref{equ::alpmulbet}, we have\begin{equation*}
		(2-2\eps)2^{2\ell}\leq e(G)\leq |\mX\cap 2^A|\cdot|\mY|+|\mX\sm 2^A|\cdot\frac{2}{3}|\mY|=\frac{2+\theta}{3}\alpha\beta 2^{2\ell}\leq \frac{2+\theta}{3}\cdot \frac{9}{4}(1+\etan)^2 2^{2\ell}.
	\end{equation*}
	%Hence, we have 
	Hence $\frac{2+\theta}{3}\cdot\frac{9}{4}(1+\etan)^2\geq 2-2\eps$, which implies that 
	$\theta\geq 2/3-3\eps$, as desired.
	Similarly, we can prove %$\phi\geq 2/3-3\eps= 2/3-o(1)$
	$\phi\geq 2/3-3\eps$.
	
	If $|A|\leq \ell-1$, then
	\begin{equation*}
		|\mX|=\frac{|\mX\cap 2^A|}{\theta}\leq \frac{|2^A|}{\theta}\leq \frac{2^{\ell-1}}{2/3-3\eps}<(1-2\eps)2^\ell,
	\end{equation*}
	a contradiction to~\eqref{eqn::bound alphabeta}, so we have $|A|\geq \ell$. Similarly, we have $|B|\geq \ell$. Therefore, $[n]=A\cup B$ is an equipartition.
\end{proof}

According to \cref{clm::size}, we can assume from now on that $|A|=\ell$ and $|B|=\ell+1$. We claim that \begin{equation}\label{eqn::bound2 alphabeta}
	\alpha\leq \frac{3}{2}+8\eps,\quad \beta\geq \frac{3}{2}-7\eps,
\end{equation} which are better bounds for $\alpha, \beta$ than~\eqref{eqn::bound alphabeta}.
Notice that we only need to show $\beta\geq 3/2-7\eps$, as $\alpha+\beta\leq 3(1+\etan)$ will then imply $\alpha\leq 3(1+\etan)-\beta\leq 3/2+8\eps$. Suppose that $\beta=3/2-\gamma$ with some $\gamma \geq 0$, then $\alpha\leq 3/2+3\etan+\gamma$ since $\alpha+\beta\leq 3(1+\etan)$ by \eqref{equ::alpplusbet}. By \eqref{eG} and \cref{obs}, we have
\begin{equation*}
	\begin{aligned}
		&(2-2\eps)2^{2\ell}\leq e(G)\leq |\mX\cap 2^A||\mY|+|\mX\sm 2^A| \cdot\frac{2}{3}|\mY|= |\mX\cap 2^A||\mY|+\left(|\mX|-|\mX\cap 2^A|\right) \cdot \frac{2}{3}|\mY|\\
		&=\left(\frac{1}{3}|\mX\cap 2^A|+\frac{2}{3}|\mX|\right)|\mY| \leq \left(\frac{1}{3}\cdot 2^\ell+\frac{2}{3}\left(\frac{3}{2}+3\etan+\gamma\right)2^\ell\right)\left(\frac{3}{2}-\gamma\right)2^\ell
		\le\left(2+3\etan-\frac{1}{3}\gamma\right)2^{2\ell},
	\end{aligned}
\end{equation*}
which implies that $2-2\eps\leq 2+3\etan-\gamma/3$. Therefore, $\gamma\leq 7\eps$.

The final %two 
three claims show that $\mX$ and $\mY$ are not too far from $2^A$ and $2^B$, respectively.
\begin{claim}\label{clm::3inES}
	(i) $|2^A\sm \mX|\leq 24\eps 2^\ell$.\\
	(ii) $|\mY\sm 2^B|\leq (\sqrt{\eps}+3\eps)2^\ell$.
\end{claim}

\begin{proof}%of Claim~\ref{clm::3inES}
	Let $$D\ce\min\{|2^B\sm \mY|, |\mY\sm 2^B|\},$$ then $D\leq |\mY\sm 2^B|\leq (1/3+3\eps)|\mY|<2|\mY|/3$ by \cref{clm::size}. Define $\mY'$ from $\mY$ by adding $D$ sets in $2^B\sm \mY$ and deleting $D$ sets in $\mY\sm 2^B$. Thus, $|\mY'|=|\mY|\leq (2+3\eps)2^\ell$ by~\eqref{eqn::bound alphabeta}.
	Note that $\mX$ and $\mY'$ are not necessarily disjoint. Let $G_1 = G_{\mX,\mY'}$.
	If $D=|2^B\sm \mY|\leq |\mY\sm 2^B|$, then $2^B\subseteq \mY'$ and $|\mY'\sm 2^B|=|\mY'|-2^{\ell+1}\leq 2\eps 2^{\ell+1}$; if $D=|\mY\sm 2^B|\leq |2^B\sm \mY|$, then $\mY'\subseteq 2^B$ and $|\mY'\sm 2^B|=0$. In both cases, we have $|\mY'\cap 2^B|\geq|\mY\cap 2^B|$ and \begin{equation}\label{eqn::size difference}
		|\mY'\sm 2^B|\leq 2\eps 2^{\ell+1}.
	\end{equation} We now compare $e(G_1)$ and $e(G)$. Every deleted $Y\in \mY\sm 2^B$ has at most $2|\mX|/3$ neighbors in $\mX$ by \cref{obs}. On the other hand, every added $S\in 2^B\sm\mY$ is disjoint from every set in $2^A\cap \mX$, thus has at least $|2^A\cap \mX|\geq (2/3-3\eps)|\mX|$ neighbors in $\mX$ by \cref{clm::size}. Therefore, by \eqref{equ::alpplusbet},
		\begin{equation*}
			e(G)-e(G_1)\leq D\left(\frac{2}{3}|\mX|-\left(\frac{2}{3}-3\eps\right)|\mX|\right)\leq \frac{2}{3} |\mY|\cdot 3 \eps|\mX|=2\eps\alpha\beta2^{2\ell}\leq 2\eps \cdot \frac{9}{4}(1+\etan)^2 2^{2\ell}\le 3\eps2^{2\ell+1},
		\end{equation*}
		which, with \eqref{eG}, implies that
		\begin{equation*}
			e(G_1)\geq e(G)-3\eps2^{2\ell+1}\geq (1-\eps)2^{2\ell+1}-3\eps2^{2\ell+1}=(1-4\eps)2^{2\ell+1}.
		\end{equation*}
		Similarly, we define another bipartite graph $G_2$. Let $$C\ce\min\{|2^A\sm \mX|, |\mX\sm 2^A|\}.$$ Define $\mX'$ from $\mX$ by adding $C$ sets in $2^A\sm \mX$ and deleting $C$ sets in $\mX\sm 2^A$. Thus, $|\mX'|=|\mX|\geq (1-3\eps)2^\ell$ by~\eqref{eqn::bound alphabeta}. Note that $\mX'$ and $\mY'$ are not necessarily disjoint.
		Let $G_2 = G_{\mX', \mY'}$.
		If $C=|2^A\sm \mX|\leq |\mX\sm 2^A|$, then $2^A\subseteq \mX'$ and $|\mX'\cap 2^A|=|2^A|=2^\ell$; if $C=|\mX\sm 2^A|\leq |2^A\sm \mX|$, then $\mX'\subseteq 2^A$ and $|\mX'\cap 2^A|=|\mX'|\geq (1-3\eps)2^\ell$. In both cases, we have $|\mX'\cap 2^A|\geq (1-3\eps)2^\ell$. We now compare $e(G_2)$ and $e(G_1)$. Every deleted $X\in \mX\sm 2^A$ intersects $B$, thus has at most $2^\ell$ neighbors in $2^B$. By~\eqref{eqn::size difference}, $X$ has at most $2\eps2^{\ell+1}$ neighbors in $\mY'\sm 2^B$, thus has at most $(1+4\eps)2^\ell$ in $\mY'$. On the other hand, every added $S\in 2^A\sm\mX$ is disjoint from every set in $2^B\cap \mY'$, thus has at least $$|2^B\cap \mY'|=|\mY'|-|\mY'\sm 2^B|\geq (3/2-7\eps)2^\ell-(2\eps)2^{\ell+1}=(3/2-11\eps)2^\ell$$ neighbors by~\eqref{eqn::bound2 alphabeta} and~\eqref{eqn::size difference}. Therefore, \begin{equation}\label{eqn::eG2l}
			e(G_2)\geq e(G_1)+C\left(\frac{3}{2}-11\eps- (1+4\eps)\right)2^\ell\geq (1-4\eps)2^{2\ell+1} + C\left(\frac{1}{2}-15\eps\right)2^\ell.
		\end{equation}
		If $|\mX'|\leq 2^\ell$, then by \eqref{equ::alpplusbet}, we have that $e(G_2)\leq |\mX'||\mY'|\leq |\mX'|(3(1+\etan)2^\ell-|\mX'|)$ attains its maximum value $(1+\frac{3}{2}\etan)2^{2\ell+1}$ when $|\mX'|=2^\ell$.
		If $|\mX'| > 2^\ell$, then let $|\mX'|=(1+\gamma)2^\ell$ for some $\gamma > 0$.
		Since $|\mX'|+|\mY'| = |\mX| + |\mY| \leq 3(1+\etan) 2^\ell$ by \eqref{equ::alpplusbet}, we have $|\mY'|\leq (2+3\etan-\gamma)2^\ell$. Recalling $|\mY'| = |\mY| \ge (3/2 - 7\eps) 2^{\ell}$ by~\eqref{eqn::bound2 alphabeta}, we have $\gamma \leq 1/2+7\eps+3\etan$.
		By the definition of $\mX'$, we have $\mX'\supseteq 2^A$, so $|\mX' \sm 2^A| = |\mX'| - |2^A| = \gamma 2^\ell$. Every $X\in \mX'\sm 2^A$ intersects $B$, thus has at most $2^\ell$ neighbors in $\mY'\cap 2^B$, so it has at most $(1+4\eps)2^\ell$ neighbors in $\mY'$ by~\eqref{eqn::size difference}.
		Therefore, \begin{equation*}
			\begin{aligned}
				e(G_2)&\leq |\mX'\sm 2^A|(1+4\eps)2^\ell+|\mX'\cap 2^A||\mY'|
				\leq \gamma 2^\ell(1+2\eps)2^\ell+2^\ell(2+3\etan-\gamma)2^\ell\\
				&=(1+\eps\gamma+\frac{3}{2}\etan)2^{2\ell+1}\leq (1+\eps)2^{2\ell+1}.
			\end{aligned}
		\end{equation*} In both cases, we have \begin{equation}\label{eqn::eG2u}
			e(G_2)\leq (1+\eps)2^{2\ell+1}. \end{equation} Combining~\eqref{eG},~\eqref{eqn::eG2l}  and~\eqref{eqn::eG2u}, we conclude that \begin{equation}\label{eqn::eG-eG1}
			e(G_1)-e(G)\leq e(G_2)-e(G)\leq (1+\eps)2^{2\ell+1}-(1-\eps)2^{2\ell+1}=2\eps2^{2\ell+1}.
		\end{equation}
	    Now we prove (i). By~\eqref{eqn::eG2l} and~\eqref{eqn::eG2u}, we have \begin{equation} \label{equ::C}
			C\leq \frac{10\eps}{1/2-15\eps}2^\ell \leq 21\eps 2^\ell,
		\end{equation}
		so if $C=|2^A\sm \mX|$, then we are done. Assume $C=|\mX\sm 2^A|$.
		Then, by~\eqref{eqn::bound alphabeta} and~\eqref{equ::C}, we have
		\begin{equation*}
			|2^A\sm \mX|=|2^A|-|2^A\cap \mX|=2^{|A|}-(|\mX|-|\mX\sm 2^A|)\leq 2^\ell-(1-3\eps)2^\ell+21\eps 2^\ell=24\eps 2^\ell,
		\end{equation*}
		as desired.
		
		For (ii), we show that it suffices to prove \begin{equation}
			D\leq \sqrt{\eps}2^\ell.
		\end{equation} Indeed, if $D=|\mY\sm 2^B|$, then we are done. Assume $D=|2^B\sm \mY|$. Then, by~\eqref{eqn::bound alphabeta}, we have
		\begin{equation*}
			|\mY\sm 2^B|=|\mY|-|2^B\cap \mY|=|\mY|-(|2^B|-|2^B\sm \mY|)\leq (2+3\eps)2^\ell-2^{2\ell+1}+\sqrt{\eps}2^\ell=\left(\sqrt{\eps}+ 3\eps \right)2^\ell,
		\end{equation*}as desired.
		
		Now suppose for a contradiction that $D> \sqrt{\eps}2^\ell$. We claim that there exists some $Y\in \mY\sm 2^B$ having at least $2|\mX|/3-8\sqrt{\eps}2^\ell$ neighbors in $\mX$. Otherwise, recall that every $S\in \mY'\sm \mY\subseteq 2^B\sm\mY$ has at least $(2/3-3\eps)|\mX|$ neighbors in $\mX$ and $|\mX|\leq (3/2+8\eps)2^\ell$ by~\eqref{eqn::bound2 alphabeta}. Then,
		\begin{equation*}
			e(G_1)-e(G)
			\geq D\left(\left(\frac{2}{3}-3 \eps \right)|\mX|-\left(\frac{2|\mX|}{3}-8\sqrt{\eps}2^\ell\right)\right)
			\geq \sqrt{\eps}\left(8\sqrt{\eps}-\frac{9}{2}\eps-24\eps^2\right)2^{2\ell}>2\eps 2^{2\ell+1},
		\end{equation*}
		a contradiction to~\eqref{eqn::eG-eG1}. Fix some $Y\in \mY\sm 2^B$ having at least $2|\mX|/3-8\sqrt{\eps}2^\ell$ neighbors in $\mX$.
		Note that there exists $i\in[n]$ such that $i\in Y\cap A$.
		We may assume without loss of generality that $i=n$. By the definition of $G$, at most $|\mX|-(2|\mX|/3-8\sqrt{\eps}2^\ell)=(1/3+o(1))|\mX|$ sets in $\mX$ contain $n$, i.e., $|\mX_n^+|/|\mX| \le 1/3+o(1)$. On the other hand, since $|2^A\sm \mX|\leq 24\eps 2^\ell$ by (i), at least $2^{\ell-1}-24\eps 2^\ell=(1-o(1))2^{\ell-1}$ subsets of $A$ containing $n$ are contained in $\mX$, so $|\mX_n^+|\ge (1-o(1))2^{\ell-1}$ and hence $|\mX| \ge (3/2-o(1))2^\ell$.
		Recalling that $|\mX| \le (3/2 + o(1))2^\ell$ by~\eqref{eqn::bound2 alphabeta}, we have
		$|\mX|=(3/2-o(1))2^\ell$, which then implies that $|\mX_n^+|=(1-o(1))2^{\ell-1}$ and $|\mY|=(3/2-o(1))2^\ell$. By the definition of $B$, we have $|\mY_n^+|/|\mY| \leq 1/3$, since $n\in A$. Recall~\eqref{eqn::derive size}, where we have now $x,y\le 1/3 + o(1)$ and $\alpha,\beta \le 3/2 + o(1)$.
		We get $y = 1/3 + o(1)$ and hence
		%Moreover, both $\mX$ and $\mY$ contain $(1-o(1))2^{\ell-1}$ sets containing $n$, i.e.,
		$|\mY_n^+|=(1-o(1))2^{\ell-1}$. Therefore, $|\mX_n^-|,|\mY_n^-|=(1-o(1))2^\ell$. Again, by \cref{lem::size}, we obtain a contradiction and complete the proof of \cref{clm::3inES}.
\end{proof}

\begin{claim} \label{cla::2by}
	$|2^B\sm \mY|\leq 6\sqrt{\eps}2^\ell$.
\end{claim}

\begin{proof}%of Claim~\ref{cla::2by}
	%Denote $2^B\sm \mY$ by $\B$ and $\mX\sm 2^A$ by $\mX_0$.
	For every $e\in E(G)$, call $e$ \emph{bad} if $e$ has an endpoint in $\mY\sm 2^B$ and \emph{good} otherwise. By \cref{clm::3inES} (ii) and~\eqref{eqn::bound2 alphabeta}, $G$ has at most
	\begin{equation} \label{equ::numBadEdges}
		|\mY\sm 2^B||\mX|\leq \left(\sqrt{\eps}+3\eps\right)2^\ell\cdot \left(3/2+8\eps\right)2^\ell\leq 2\sqrt{\eps}2^{2\ell}
	\end{equation}
	bad edges.
	
	Fix $S'\in 2^B\sm \mY$ and choose a set $F\subseteq [n]$ of the form $F=S\cup S'$, where $S\subseteq A$.
	%Since $\mH$ generates all but at most $\eps2^n$ subsets of $[n]$, at least $2^\ell-\eps 2^n$ subsets of the form $S\cup S'\subseteq [n]$, where $S\subseteq A$, can be expressed as a union of two disjoint sets of $\mH$.
	If $F$ corresponds to a good edge of $G$, assume that $F = X\cup Y$ with $X\in\mX, Y\in\mY \cap 2^B, X \cap Y=\emptyset$.
	By the assumption that $S' \in 2^B \sm \mY$ and $Y \in \mY \cap 2^B$, we have $Y \neq S'$ and hence $Y \subsetneq S'$. Hence, $X \cap B \neq \emptyset$, so $X \in \mX \sm 2^A$. Furthermore, $X \cap A = S$,
	%Hence, $Y\subsetneq S'$ is a proper subset of $B$ since $S'\notin\mY$, which implies that $X\in \mX\sm 2^A$ since $X$ must intersect $B$. Therefore, we have $X\cap A=S$,
	which implies that $S$ is determined by $X$.
	%a different $S$ would correspond to a different $X$. There are at most
	By \eqref{eqn::bound2 alphabeta} and \cref{clm::3inES}~(i), we have
	$$
	|\mX\sm 2^A|=|\mX|-|\mX\cap 2^A|= |\mX| - (|2^A| - |2^A \sm \mX|) \leq (3/2+8\eps)2^\ell-(1-24\eps) 2^\ell=(1/2+32\eps)2^\ell.
	$$
	%different $X$'s
	Hence, for every fixed $S'\in 2^B\sm \mY$, there are at most $(1/2+32\eps)2^\ell$ sets $F$ of the form $F=S\cup S'$, where $S\subseteq A$, corresponding to good edges of $G$. There are $2^\ell|2^B\sm \mY|$ sets of the form $S\cup S'$ with $S\subseteq A, S'\in 2^B\sm \mY$ in total, so there are at least $|2^B\sm \mY|(2^\ell-(1/2+32\eps)2^\ell)$ sets of the form $S\cup S'$ with $S\subseteq A, S'\in 2^B\sm \mY$ that correspond to bad edges of $G$ or do not correspond to edges of $G$.
	Recalling that $\mH$ is a $(1-\eps)$-$2$-generator for $[n]$ and using~\eqref{equ::numBadEdges}, we have
	\begin{equation*}
		|2^B\sm \mY|\left(2^\ell-\left(1/2+32\eps\right)2^\ell\right)\leq 2\sqrt{\eps}2^{2\ell}+\eps 2^n=\left(2\sqrt{\eps}+2\eps\right)2^{2\ell}.
	\end{equation*}
	Therefore, we get $|2^B\sm \mY|\leq 6\sqrt{\eps}2^\ell$.
\end{proof}

\begin{claim} \label{cla::x2a}
	$|\mX \sm 2^A| \le 7\sqrt{\eps} 2^\ell$.
\end{claim}
\begin{proof}%of Claim~\ref{cla::x2a}
	By Claims~\ref{clm::3inES}~(i) and~\ref{cla::2by}, we have
	\begin{align*}
		3(1+\etan)2^\ell \ge |\mH|
		&=  |\mX| +  |\mY|
		= |\mX \sm 2^A| + \left( |2^A| - |2^A \sm \mX| \right) + |\mY \sm 2^B| + \left( |2^B| - |2^B \sm \mY| \right) \\
		&\ge |\mX \sm 2^A| + (1-24\eps)2^\ell + (2 - 6\sqrt{\eps})2^\ell = |\mX \sm 2^A| +3\cdot 2^\ell -(6\sqrt{\eps} + 24\eps)2^\ell,
	\end{align*}
	implying $|\mX \sm 2^A| \le 7\sqrt{\eps} 2^\ell$.
\end{proof}

By \cref{clm::size}, $[n] = A \cup B$ is an equipartition. By Claims~\ref{clm::3inES},~\ref{cla::2by} and~\ref{cla::x2a}, we have $$|\mH\Delta(2^A \cup 2^B)| \le (14\sqrt{\eps} + 27\eps)2^\ell \le \eps' 2^\ell,$$ completing the proof of \cref{thm::modify}.
\end{proof}

\section{Proof of \cref{thm::main}} \label{sec::pro}
In this section, we will employ a perturbation argument to deduce \cref{thm::main} from \cref{mlem}.
Let $n$ be a sufficiently large integer and $\mF$ be a maximal $3$-wise intersecting family on $[n]$ of size at most $2^{\lfloor n/2\rfloor}+2^{\lceil n/2 \rceil}-3$. Let $\mH:=\bar{\mF}=\{F^c: F\in\mF\}$ and fix some $\eps\in(0,1/4)$. By \cref{mlem}, there exists $S\subseteq [n]$ of size $\lfloor n/2 \rfloor$ such that $\mF_0\ce\{A: A\subseteq S\}\cup \{B: B\subseteq S^c\}$ satisfies $|\mH\Delta\mF_0|\leq \eps 2^{\lfloor n/2 \rfloor}$.

Recall that $\mH$ is downward closed since $\mF$ is upward closed, so $\emptyset\in\mH$. We first prove that $S\notin \mH=\bar{\mF}$. Suppose for a contradiction that $S^c\in \mF$. Since $|\mH\Delta\mF_0|\leq \eps 2^{\lfloor n/2 \rfloor}$, among the $2^{\lceil n/2 \rceil}$ subsets of $S^c$, there exists some $A\subseteq S^c$ such that both $A$ and $S^c\sm A$ are contained in $\mH$. However, this would imply that $S^c, A^c, (S^c\sm A)^c=S\cup A$ are in $\mF$. As $S^c\cap A^c\cap (S\cup A)=\emptyset$, this contradicts that $\mF$ is a $3$-wise intersecting family. It can be proved similarly that $S^c\notin \mH$.

We work with the following partition $\mH=\mH_1\cup\mH_2\cup\mH_3\cup\{\emptyset\}$, where 
\[
\mH_1\ce\{A\in\mH\colon \emptyset\subsetneq A\subsetneq S\}, \quad \mH_2\ce\{B\in\mH\colon \emptyset\subsetneq B\subsetneq S^c\}, \quad \mH_3\ce\mH\sm\mF_0.
\]

\begin{claim}\label{clm::mpf1}
	$|2^{[n]}\sm (\mF\cup \mF_0)|\leq |\mH_1||\mH_2|+|\mH_3|\cdot\eps 2^{\lfloor n/2 \rfloor}$.
\end{claim}

\begin{proof}%of \cref{clm::mpf1}
	If $A\notin \mF$, then as stated in the proof of \cref{lem::toGen}, there exist $B,C\in \mH$ such that $A=B\cup C$ with $B\cap C=\emptyset$. %Among those, 
 There exists a pair $\{B,C\}$ such that $$g(\{B,C\})\ce\min \{|B\sm S|+|C\sm S^c|, |B\sm S^c|+|C\sm S|\}$$ attains its minimum value over all such choices of $B$ and $C$. Therefore, we can define an injection $f$ on $2^{[n]}\sm (\mF\cup 2^S\cup 2^{S^c})$ by mapping $A$ to a pair of sets $\{B,C\}$ in $\mH$, which has minimum $g(\{B,C\})$ given $A=B\cup C$ and $B\cap C=\emptyset$. Since $\emptyset\in \mH$, $S,S^c\notin\mH$ and $\mH$ is %closed downwards
downward closed, we have that $f(A)=\{B,C\}$ must be one of the following two types: $|\{B,C\}\cap \mH_1|=|\{B,C\}\cap \mH_2|=1$; $|\{B,C\}\cap \mH_3|\geq 1$. The number of pairs $\{B,C\}$ satisfying $|\{B,C\}\cap \mH_1|=|\{B,C\}\cap \mH_2|=1$ is at most $|\mH_1||\mH_2|$, so it suffices to prove that if $f(A)=\{B,C\}$ with $B\in \mH_3$, then there are at most $\eps 2^{\lfloor n/2 \rfloor}$ choices for $C$.
	
	Suppose $B=X\cup Y$, where $\emptyset\neq X\subsetneq S, \emptyset\neq Y\subsetneq S^c$. Then, $g(\{B,C\})>0$ and $Y\in \mH$, since $Y\subseteq B\in\mH$ and $\mH$ is downward closed. There are three possibilities for $C$.
	
	(1) $C\subseteq S$: Define $B' \ce X\cup C, C'\ce Y$, then $B'\cup C'=B\cup C=A$. Since $g(\{B',C'\})=0<g(\{B,C\})$ and $C'=Y\in\mH$, we have that $B'\notin \mH$ by the definition of $f$. Note that $B'=X\cup C$ is determined by $C$, so the number of choices for $C\subseteq S$ is at most $|2^S\sm \mH|$.
	
	(2) $C\subseteq S^c$: Similarly, the number of choices for $C$ is at most $|2^{S^c}\sm\mH|$.
	
	(3) $C\in \mH\sm\mF_0$: The number of choices for such $C$ is at most $|\mH\sm\mF_0|$.
	
	In summary, the number of choices for $C$ is at most \begin{equation*}
		|2^S\sm \mH|+|2^{S^c}\sm \mH|+|\mH\sm\mF_0|=|\mH\Delta\mF_0|\leq \eps 2^{\lfloor n/2 \rfloor},
	\end{equation*} as desired.
\end{proof}

\begin{claim}\label{clm::mpf2}
	$|\mH_1||\mH_2|+|\mH_3|\cdot\eps 2^{\lfloor n/2 \rfloor}
	\leq (2^{\lfloor n/2 \rfloor}-2)(2^{\lceil n/2 \rceil}-2)$. Equality holds if and only if $\mH_1=\{A\colon \emptyset\subsetneq A\subsetneq S\}$, $\mH_2=\{B\colon \emptyset\subsetneq B\subsetneq S^c\}$ and $\mH_3=\emptyset$.
\end{claim}

\begin{proof}%of \cref{clm::mpf2}
	By the definitions of $\mH_1$ and $\mH_2$, we have
	\begin{equation}\label{eqn::g1g2}
		|\mH_1||\mH_2|\leq (2^{\lfloor n/2 \rfloor}-2)(2^{\lceil n/2 \rceil}-2).
	\end{equation}
	Notice that $\min\{|\mH_1|,|\mH_2|\}\geq 2^{\lfloor n/2 \rfloor}-|\mF_0\sm \mH|\geq 2^{\lfloor n/2 \rfloor}-|\mF_0\Delta\mH|\geq (1-\eps)2^{\lfloor n/2 \rfloor}$.
	Define the function $\varphi(h_1,h_2,h_3)\ce h_1h_2+h_3\cdot\eps2^{\lfloor n/2 \rfloor}$. Given that $\min\{h_1,h_2\}\geq(1-\eps)2^{\lfloor n/2 \rfloor}$, we claim that $\varphi(h_1',h_2,h_3')> \varphi(h_1,h_2,h_3)$ if $h_1'=h_1+1$ and $h_3'=h_3-1$. %$h_1'+h_2'=h_1+h_2+1$ and $h_3'=h_3-1$.
	In fact, we have% since $h_1'h_2'\in\{(h_1+1)h_2,h_1(h_2+1)\}$, we have}
	\[
	\varphi(h_1',h_2,h_3')-\varphi(h_1,h_2,h_3)= h_2-\eps 2^{\lfloor n/2 \rfloor}\geq (1-2\eps)2^{\lfloor n/2 \rfloor}>0.
	\] Therefore, fixing $h_1+h_2+h_3$ and stipulating $h_1,h_2\geq (1-\eps)2^{\lfloor n/2 \rfloor}$ and $h_3\geq 0$, the function $\varphi(h_1,h_2,h_3)$ attains its maximum when $h_3=0$. Take $h_i=|\mH_i|$ for $i\in[3]$. Combining with~\eqref{eqn::g1g2}, we have \[
	|\mH_1||\mH_2|+|\mH_3|\cdot\eps 2^{\lfloor n/2 \rfloor}\leq (2^{\lfloor n/2 \rfloor}-2)(2^{\lceil n/2 \rceil}-2),\] where equality holds if and only if $\mH_1=\{A\colon \emptyset\subsetneq A\subsetneq S\}$, $\mH_2=\{B\colon \emptyset\subsetneq B\subsetneq S^c\}$ and $\mH_3=\emptyset$.
	%By \crefeq{clm::mpf1}, it suffices to prove that as long as there exists some $A\in \mH_3$ with $A\notin 2^S\cup 2^{S^c}$, we can replace $A$ by some $B\in \F_0$. In fact, since $|\mH|\leq 3\cdot2^\ell-3$, there exists some $B\in \F_0$ such that $B\notin \mH$. We replace $A$ by $B$ in $\mH$ and consider the change of $|\mH_1||\mH_2|+|\mH_3|\cdot\eps 2^\ell$. Since deleting $A$ from $\mH$ would decrease $|\mH_1||\mH_2|+|\mH_3|\cdot\eps 2^\ell$ by at most $\eps 2^\ell$, and adding $B$ to $\mH$ would increase $|\mH_1||\mH_2|+|\mH_3|\cdot\eps 2^\ell$ by at least $2^{\ell+1}-|\F_0\sm\mH|\geq (1-\eps)2^\ell$, we are done.
\end{proof}

By Claims~\ref{clm::mpf1} and~\ref{clm::mpf2}, we have
\begin{equation*}
	2^{n}-(2^{\lfloor n/2 \rfloor} + 2^{{\lceil n/2 \rceil}}-1)-|\mF|
	\leq |2^{[n]}\sm (\mF\cup \mF_0)|
	\leq (2^{\lfloor n/2 \rfloor}-2)(2^{\lceil n/2 \rceil}-2)
	=2^{n}-2\cdot 2^{\lfloor n/2 \rfloor} - 2 \cdot 2^{\lceil n/2 \rceil}+4,
\end{equation*}
which implies that $|\mF|\geq 2^{\lfloor n/2 \rfloor} + 2^{\lceil n/2 \rceil}-3$. By our assumption that $|\mF|\leq 2^{\lfloor n/2 \rfloor} + 2^{\lceil n/2 \rceil}-3$, equality has to hold in \cref{clm::mpf2}. Hence, we have $\mH=\mF_0$, which means that $\mF$ is a balanced pair of linked cubes. This completes our proof of \cref{thm::main}.

\section{The Case when $k\geq 4$}\label{sec::largek}

Denote by $f(n,k)$ the minimum possible size of a maximal $k$-wise intersecting family on $[n]$. For the case when $k\geq 4$, we remark that our method also gives the following general bounds for $f(n,k)$.

\begin{proposition}\label{prop: larger k}
	For every $k\geq 4$  there exist positive constants $c_k$ and $d_k$ such that for every positive integer $n$, we have \[
	c_k\cdot2^{n/(k-1)}\leq f(n,k)\leq d_k\cdot2^{n/\lceil k/2\rceil}.\]
\end{proposition}

\begin{proof}
	The proof of \cref{lem::toGen} easily generalizes to give a lower bound for $f(n,k)$. %the size of a maximal $k$-wise intersecting family.
	%In particular,
	Indeed, there is an injective map from $\mF^c$ to $\binom{\mF}{k-1}$, such that each set $A \in \mF^c$ is sent to a $(k-1)$-tuple $\{B_1, B_2, \ldots, B_{k-1}\}$ of sets in $\mF$ such that $B_1 \cap B_2 \cap \cdots \cap B_{k-1} = A^c$ (unless $A$ is so small that fewer than $k-1$ sets contain $A^c$). This observation immediately leads to the lower bound.
	
	The linked cubes construction also generalizes.
	In more detail, suppose that $\ell$ divides $n$, and let $S_1, S_2, \ldots, S_\ell$ be a partition of $[n]$ such that $|S_i| = n/\ell$ for each $i$.
	Let $\mF_i = \{A\colon [n]\sm S_i \subsetneq A \subseteq [n]\}$, and let $\mF = \mF_1 \cup \mF_2 \cup \cdots \cup \mF_\ell$.
	For any set $\mG$ of $2\ell-1$ elements of $\mF$, there is an $i$ such that $|\mG \cap \mF_i| \leq 1$. The sets in $\mG-\mF_i$ all contain $S_i$, and the set in $\mG \cap \mF_i$ contains at least one element of $S_i$, which element is in every set of $\mG$. Hence, $\mG$  is $(2\ell-1)$-wise intersecting.
	On the other hand, if $A \notin \mF$, then either $A$ is missing elements from two distinct sets $S_i$, or $A$ is missing all of the elements from some $S_i$.
	In either case, it is easy to find a  non-intersecting $(2\ell - 1)$-tuple of elements in $\mF \cup \{A\}$.
	
	For the case when $k$ is even, let $\mF'$ be a maximal $(k-1)$-intersecting family on $[n-1]$.
	Then the family $\mF = \{A \cup \{n\}: A \in \mF\} \cup \{[n-1]\}$ is a maximal $k$-wise intersecting family.
	
	These constructions give the upper bound.
\end{proof}

Very recently, Janzer~\cite{Jan21} showed that the lower bound in \cref{prop: larger k} is the correct order of magnitude of $f(n,k)$, by constructing a maximal $k$-wise intersecting family of size $O(2^{n/(k-1)})$ for every $k\geq 3$ and $n$. Note that in the special case $k=3$, \cref{thm::main} matches Janzer's~\cite{Jan21} construction.

Assume that $n$ is sufficiently large. Let $\mF$ be a maximal $k$-wise intersecting family on $[n]$ with minimum size. Similarly to the case $k=3$, one can show that $\bar{\mF}$ is a $(1-\eps)$-$(k-1)$-generator for $[n]$, where $\eps=o(1)$. %\cref{thm::stability n even} and 
A modification of the method in %\cref{sec::pro}
this paper could be used to determine the structure of $\bar{\mF}$, if $|\bar{\mF} \sm \mF_{n,k-1}|$ was small.
%its size was about the same as $|\mF_{n,k}|$.
However, Janzer~\cite{Jan21} showed that %it is not the case.
$|\bar{\mF} \sm \mF_{n,k-1}|$ cannot be small.

\begin{theorem}[{\cite[Lemma 1.2]{Jan21}}]\label{janzer1.2}
	%Let $k\geq 4$ and $\mF_0$ be a series of $k$ cubes on $[n]$, i.e., there exists a partition $S_1\cup\ldots\cup S_{k-1}$ of $[n]$ with $|S_i|=n/(k-1)+O(1)$ for each $i\in[k-1]$ and $\mF_0=2^{S_1}\cup\ldots\cup 2^{S_{k-1}}$. If $|\mF\sm \mF_0|\leq c\cdot 2^{n/(k-1)}$ and $n$ is large enough compared with $c$, then $\bar{\mF}$ cannot be maximal $k$-wise intersecting.
	For every $k \ge 4$, there exist $d= d(k) > 0$, $c = c(k) >0$ and $n_0 = n_0(k) >0$ such that the following holds when $n \ge n_0$. Let $S_1 \cup \cdots \cup S_{k-1}$ be a partition of $[n]$, where $\frac{n}{k-1} -d \le |S_i| \le \frac{n}{k-1} +d$ for every $i \in [k-1]$. Let $\mF_0=2^{S_1}\cup\cdots\cup 2^{S_{k-1}}$. For every set system $\mF \subseteq 2^{[n]}$ with $|\mF \sm \mF_0| \le c\cdot 2^{n/(k-1)}$, $\bar{\mF}$ cannot be maximal $k$-wise intersecting.
\end{theorem}

Combining our method with %\ref{thm::stability n even} and 
\cref{janzer1.2}, we obtain the following result.

\begin{proposition}
	For every $k \ge 4$, there exists $c = c(k)>0$ such that
	$$
	f(n,k)\geq (1+c)|\mF_{n,k-1}| = (1+c) \left((k-1)2^{\frac{n}{k-1}} -k + 2\right),
	$$
	when $n$ is divisible by $k-1$.
\end{proposition}
Therefore, a maximal $k$-wise intersecting family of minimum size will necessarily have a more complex structure when $k \ge 4$.
It is worth mentioning that the exact value of the upper bound on Janzer's construction~\cite{Jan21} is $(k-1)2^{k-3}2^{n/(k-1)}-(k-2)(2^{k-1}-1)$, which is %larger than $|\mF_{n,k-1}|$ with about a multiplicative factor of $2^{k-3}$.
about $2^{k-3}|\mF_{n,k-1}|$.

\section*{Acknowledgment}

The authors are grateful for Jingwei Xu, Simon Piga and Andrew Treglown, who participated in fruitful discussions at the  beginning of  the project. Simon Piga and Andrew Treglown's visit to University of Illinois was partially supported by NSF RTG grant DMS 1937241.

%\section{Discussions}

%When $n$ is a multiple of $k$, we have the stability but...

%When $n$ is not a multiple of $k$, it is hard to deduce a similar lemma as \cref{lem::bipartite}, as ...

\end{document}